\input amstex

\documentstyle{amsppt}
\define\mL{L\kern-0.08cm\char39}
\magnification=\magstep1
\topmatter
\title
The graph topology
\endtitle
\rightheadtext{}
\author  \v Lubica Hol\'a
\endauthor
\abstract We study topological properties of the graph topology on the space of continuous real-valued functions.

\endabstract

\bigskip
\bigskip
\bigskip

\address Academy of Sciences, Institute of Mathematics, \v
Stef\'anikova 49,
81473 Bratislava, Slovakia
\endaddress
\email hola\@mat.savba.sk
\endemail

\keywords continuous real-valued function, Vietoris topology, graph topology.
\v L. Hol\'a would like to thank to grant Vega 2/0047/10
\endkeywords

\endtopmatter

\document
\bigskip
\bigskip
\bigskip
\bigskip
\bigskip
\heading 1. Preliminaries
\endheading
\bigskip
\bigskip
Let $X$ be a topological space and $C(X)$ be the space of continuous real-valued functions. The graph topology $\tau_\Gamma$ was introduced by Naimpally [Na] and has as its basic open sets, sets of the form $\{f \in C(X): f \subset G\}$, where $G$ is an open subset of $X \times R$ (here $f$ is identified with its graph). If $G$ is an open set in $X \times R$ denote by $F_G$ the set $\{f \in C(X): f \subset G\}$. It is known that for a  $T_1$ space $X$ the graph topology on $C(X)$ coincides with the Vietoris topology on $C(X)$.

\bigskip
\bigskip
Najprv by sme asi mali uvazovat graph topology na vseobecnych $X, Y$.

\bigskip

Let $X$ and $Y$ be topological spaces. Denote the graph topology on $C(X,Y)$, the space of continuous functions from $X$ to $Y$,  by $\tau_{\Gamma}$.

If $X$ is a $T_1$ space then the Vietoris topology and the graph topology on $C(X,Y)$ coincide:

(For every $x \in X$ and $U$ open in $Y$, let $W(x,U) = \{f \in C(X,Y): f(x) \in U\}.$ Since $X$ is $T_1$, the set $V = (X \times U) \cup (X \setminus \{x\}) \times Y$ is open in $X \times Y$ and clearly $F_V = W(x,U)$. Now let $G \times U$ be a product of open sets in $X$ and $Y$, respectively and consider the set $(G \times U)^-$. It is easy to verify that

\centerline{$(G \times U)^- = \cup_{x \in G} W(x,U)$.}

\bigskip
The separation axioms $T_1$ and $T_2$ on $(C(X,Y),\tau_\Gamma)$ are solved in [Na]. In what follows let $X, Y$be Hausdorff spaces. It is easy to verify that if $X \times Y$ is a normal space, then $(C(X,Y),\tau_{\Gamma})$ is regular. Let $G$ be an open set in $X \times Y$ and $f \in C(X,Y)$ such that $f \in F_G$. Since the graph of $f$ is a closed set in $X \times Y$ and $f \subset G$, the normality of $X \times Y$ implies that there is an open set $H$ in $X \times Y$ such that

\centerline{$f \subset H \subset \overline H \subset G$.}

\bigskip
Since the set $F_{\overline H} = \{h \in C(X,Y): h \subset \overline H\}$ is closed in $(C(X,Y),\tau_\Gamma)$, we are done.

\bigskip
The condition of normality of $X \times Y$ is essential.

\bigskip
\proclaim{Example} Let $W_0$ be the ordinals less than or equal to the first uncountable ordinal $\Omega$ with the order topology and $W = W_0 \setminus \{\Omega\}$, with the induced topology from $W_0$. It is known that $W \times W_0$ is not normal. We will show that $(C(W,W_0),\tau_{\Gamma})$ is not regular. Let $f: W \to W_0$ be the identity function. The set $W \times W$ is open in $W \times W_0$ and $f \subset W \times W$. There is no open set $G$ in  $W\times W$ such that

\centerline{$f \subset G \subset \overline G \subset W \times W$.}
\bigskip
It is easy to verify that for every open set $G \subset W \times W$ with $f \in F_G$ we have that the closure of $F_G$ in the graph topology is not contained in $F_{W \times W}$. Let $G$ be an open set in $W \times W_0$ such that

\centerline{$f \subset G \subset W \times W$.}

\bigskip

There is $\alpha \in W$ such that $(y,z) \in G$,  if $y > \alpha$, $z > \alpha$. Define a function $g$ as follows: $g(x) = x$ for $x \le \alpha + 1$ and $g(x) = \Omega$ for every $x > \alpha + 1$. It is easy to verify that the function $g$ is in the closure of $F_G$ in the graph topology and of course $g \notin F_{W \times W}$.

Thus $(C(W,W_0),\tau_{\Gamma})$ is not regular.

\endproclaim

\bigskip

Denote further by $LSC(X,(0,1))$ the space of all lower semicontinuous functions defined on $X$ with values in the open interval $(0,1)$.

\bigskip
\bigskip
\proclaim{Lemma 1.1} Let $X$ be a topological space. $(C(X),\tau_\Gamma)$ has as a base all sets of the form
\bigskip
\centerline{$B(f,\epsilon) = \{g \in C(X): \mid f(x) - g(x) \mid < \epsilon(x)\}$,}
\bigskip
where $f \in C(X)$ and $\epsilon \in LSC(X,(0,1))$.

\endproclaim
\demo{Proof} Let $G$ be an open set in $X \times R$ with $f \subset G$. It is easy to verify that there is an open set $U$ in $X \times R$ such that $f \subset U \subset G$, $U$ is locally bounded and $U(x)$ is an open interval for every $x \in X$. For every $x \in X$ put
\bigskip
\centerline{$\epsilon(x) = min \{\mid f(x) - y\mid: y \in U(x)^c\}$.}
\bigskip
It is easy to verify that $\epsilon$ is a positive lower semicontinuous function. Put $\eta(x) = min \{\epsilon(x), 1/2\}$. Then $\eta \in LSC(X,(0,1))$ and of course $B(f,\eta) \subset U \subset G$.
\enddemo

\bigskip
\proclaim{Lemma 1.1'} Let $X$ be a topological space and $(Y,d)$ be a metric space. Then $(C(X,Y),\tau_{\Gamma})$ has as a base all sets of the form 

\bigskip
\centerline{$B(f,\epsilon) = \{g \in C(X,Y): d(f(x),g(x)) < \epsilon(x), \forall x \in X\},$}

\bigskip

where $f \in C(X,Y)$ and $\epsilon \in LSC(X,(0,1))$.
\endproclaim
\demo{Proof} Let $V$ be an open set in $X \times Y$ with $f \subset G$. For every $x \in X$ put
\bigskip

\centerline{$V(x) = \{y \in Y: (x,y) \in V\}.$}

Without loss of generality we can suppose that $V(x)^c$ is nonempty for every $x \in X$. Now define the functions $\sigma, \eta$ as follows:

\bigskip

\centerline{$\sigma(x) = d(f(x),V(x)^c)$ for every $x \in X$, and}

\bigskip

\centerline{$\eta(x) = sup \{inf \{\sigma(z): z \in U\}: U \in \Cal U(x)\}$ for every $x \in X$,}

\bigskip

where $\Cal U(x)$ is the neighborhood base at $x$. Put $\epsilon(x) = min \{\eta(x), 1\}$ for every $x \in X$. It is easy to verify that $\epsilon$ is lower semicontinuous and positive and $B(f,\epsilon) \subset F_V$.

We need to prove that $B(f,\epsilon) \in \tau_{\Gamma}$ for $f \in C(X,Y)$ and $\epsilon \in LSC(X,(0,1))$. We prove that the set

\bigskip
\centerline{$H(f,\epsilon) = \{(x,t) \in X \times Y: d(f(x),t) < \epsilon(x), x \in X\}$}

\bigskip
is open in $X \times Y$. Let $(x,t) \in H(f,\epsilon)$. Put $\alpha(x) = \epsilon(x) - d(f(x),t)$. There is an open  neighborhood $O_x$ of $x$ such that $d(f(x),f(z)) < \alpha(x)/3$ and $\epsilon(z) > \epsilon(x) - \alpha(x)/3$ for every $z \in O_x$. We claim that 

\bigskip
\centerline{$H = O_x \times \{y \in Y: d(y,t) < \alpha(x)/3\} \subset H(f,\epsilon).$}

\bigskip

Let $(z,y) \in H$. Then $d(f(z),y) \le d(f(z),f(x)) + d(f(x),y) \le d(f(z),f(x)) + d(f(x),t) + d(t,y) \le \alpha(x)/3 + d(f(x),t) + \alpha(x)/3 = 2\alpha(x)/3 + d(f(x),t) = 2/3(\epsilon(x) - d(f(x),t)) + d(f(x),t) = 2\epsilon(x)/3 + d(f(x),t)/3 < \epsilon(z)$. Thus $d(f(z),y) < \epsilon(z)$, so  we are done.

\enddemo

\bigskip

\proclaim{Proposition 1.2} Let $X$ be a topological space. The following are equivalent:

(1) $\tau_\omega = \tau_\Gamma$ on $C(X)$;

(2) for every $\epsilon \in LSC(X,(0,1))$ there is $\eta \in C(X,(0,1))$ such that $\eta(x) \le \epsilon(x)$ for every $x \in X$.
\endproclaim
\bigskip
Of course, if $X$ is countably compact or countably paracompact normal space, then the property (2) in Proposition 1.2 is satisfied.
It would be interesting to know whether the property (2) characterizes countable paracompactness in the class of normal spaces.
\bigskip
\bigskip
If $X$ is a $T_1$-space, then also $(C(X),\tau_\Gamma)$ is a $T_1$-space. By Lemma 1.1 we know that $(C(X),\tau_\Gamma)$ is a topological group.
Thus if $X$ is a $T_1$-space, $(C(X),\tau_\Gamma)$ is a Tychonoff space. $(C(X),\tau_\Gamma)$ is a topological vector space if and only if $X$ is countably compact.

\bigskip
\bigskip
The family $\{B(\epsilon): \epsilon \in LSC(X,(0,1))\}$ of sets of the form $B(\epsilon) = \{(f,g) \in C(X) \times C(X); \mid f(x) - g(x) \mid < \epsilon(x), \forall x \in X\}$ clearly satisfies the conditions for a base for a uniformity on $C(X)$ ([Ke]). Since obviously, $B(\epsilon)[f] = B(f,\epsilon)$, the uniformity induces the graph topology $\tau_\Gamma$ on $C(X)$. This uniformity will be called the uniformity of the graph topology and will be denoted by $\Cal U_\Gamma$.
\bigskip
\bigskip
\proclaim{Proposition 1.3} Let $X$ be a topological space. $(C(X),\Cal U_\Gamma)$ is a complete uniform space.
\endproclaim
\demo{Proof} Let $\{f_\lambda: \lambda \in \Lambda\}$ be a $\Cal U_\Gamma$-Cauchy net in $C(X)$. Since $\Cal U_\Gamma$ contains the uniformity $\Cal U_u$ of uniform convergence, the net is also $\Cal U_u$-Cauchy and therefore is uniformly convergent to some $f \in C(X)$. The rest of the proof is similar to the proof of an analogous statement for $\Cal U_u$. For $\epsilon \in LSC(X,(0,1))$ there is $\lambda_0 \in \Lambda$ such that for each $\lambda_1 \ge \lambda_0$, each $\lambda_2 \ge \lambda_0$ and each $x \in X$ we have
\bigskip
\centerline{$\mid f_{\lambda_1}(x) - f_{\lambda_2}(x) \mid < \epsilon(x)/2.$}
\bigskip
Fix $\lambda \ge \lambda_0$ and $x \in X$. There is $\eta \ge \lambda$ such that $\mid f_\eta(x) - f(x)\mid < \epsilon(x)/2$. Thus we have
\bigskip
\centerline{$\mid f_\lambda(x) - f(x) \mid \le \mid f_\lambda(x) - f_\eta(x) \mid + \mid f_\eta(x) - f(x) \mid < \epsilon(x).$}
\bigskip
Consequently $f_\lambda \in B(f,\epsilon)$ for all $\lambda \ge \lambda_0$, so the net $\{f_\lambda: \lambda \in \Lambda\}$ $\tau_\Gamma$-converges to $f$, which completes the proof.

\enddemo
\bigskip
\bigskip
\bigskip
\heading   2. Topological properties of the graph topology
\endheading
\bigskip
\bigskip

\proclaim{Proposition 2.1} Let $X$ be a Tychonoff space. The following are equivalent:

(1) $(C(X),\tau_\Gamma)$ is first countable;

(2) $(C(X),\tau_\Gamma)$ is metrizable;

(3) $(C(X),\tau_\Gamma)$ is completely metrizable;

(4) $(C(X),\tau_\Gamma)$ is \v Cech-complete;

(5) $X$ is countably compact;

(6) $(C(X),\tau_\Gamma)$ is Frechet;

(7) $(C(X),\tau_\Gamma)$ has a countable tightness.

\endproclaim
\demo{Proof} Since the equivalence $(1) \Leftrightarrow (2) \Leftrightarrow (3) \Leftrightarrow (4) \Leftrightarrow (5)$ is proved in [HH],
it is sufficient to prove that $(7) \Rightarrow (5)$.
Suppose $X$ is not countably compact.  There is an infinite set $\{x_n: n \in \omega\}$ in $X$ without an accumulation point. It is easy to verify that there is a sequence of open neighborhoods $\{O(x_n): n \in \omega\}$ of the points $\{x_n: n \in \omega\}$ such that the sequence $\{\overline{O(x_n)}: n \in \omega\}$ is pairwise disjoint.

Let $f_0$ be the function identically equal to $0$. Then $f_0$ belongs to the $\tau_\Gamma$-closure of the set
\bigskip
\centerline{$L = \{g \in C(X): g(x_n) > 0, \forall n \in \omega\}.$}
\bigskip
Indeed, let $\epsilon \in LSC(X,(0,1)$ and consider the set $B(f_0,\epsilon)$. Put
\bigskip
\centerline{$V(f_0,\epsilon) = \{(x,t): - \epsilon(x) < t < \epsilon(x)\}.$}
\bigskip

Then $V(f_0,\epsilon)$ is an open set in $X \times R$. For every $n \in \omega$, there is an open set $U(x_n) \subset O(x_n)$ in $X$ with $x_n \in U(x_n)$ and an open interval $(a_n,b_n)$ containing $0$ with
\bigskip

\centerline{$U(x_n) \times (a_n,b_n) \subseteq V(f_0,\epsilon)$.}

\bigskip

For every $n \in \omega$, let $H(x_n)$ be an open set in $X$ such that $x_n \in H(x_n) \subset \overline{H(x_n)} \subset U(x_n)$ and let $h_n: \overline{H(x_n)} \to [0,b_n/n]$ be a continuous function such that $h_n(x_n) = b_n/n$ and $h_n(z) = 0$ for every $z \in \overline{H(x_n)} \setminus H(x_n)$. Now let $h: X \to [0,1]$ be defined as follows: $h(x) = h_n(x)$ for every $x \in \overline{H(x_n)}$ and $h(x) = 0$ otherwise. It is easy to verify that $h$ is continuous function and $h \in B(f_0,\epsilon)$.

The countable tightness of $(C(X),\tau_\Gamma)$ implies that there is a countable subfamily $L' = \{f_n: n \in \omega\}$ of $L$ such that $f_0$ is in the $\tau_\Gamma$-closure of $L'$. Let $\eta$ be a positive function defined on $X$ as follows: $\eta(x_n) = f_n(x_n)/(n+1)$ for every $n \in \omega$ and $\eta(x) = 1/2$ otherwise. It is easy to verify that $\eta \in LSC(X),(0,1))$.

Since $f_n \notin B(f_0,\eta)$ for every $n \in \omega$, we have a contradiction.

\enddemo

\bigskip
Myslim, ze plati aj nasledujuce:
\bigskip
\proclaim{Proposition 2.1'} Let $X$ be a Tychonoff space and $Y$ be a metrizable space which contains an arc. The following are equivalent:

(1) $(C(X,Y),\tau_{\Gamma})$ is metrizable;

(2) $(C(X,Y),\tau_{\Gamma})$ is first countable;

(3) $(C(X,Y),\tau_{\Gamma})$ is Frechet;

(4) $(C(X,Y),\tau_{\Gamma})$ has a countable tightness;

(5) $X$ is countably compact.

\endproclaim

\bigskip
\proclaim{Proposition 2.2} Let $X$ be a Tychonoff space. The following are equivalent:

(1)  $(C(X),\tau_\Gamma)$ is second countable;

(2)  $(C(X),\tau_\Gamma)$ has a countable network;

(3)  $(C(X),\tau_\Gamma)$ is separable;

(4)  $(C(X),\tau_\Gamma)$ has a countable chain condition;

(5)  $X$ is compact and metrizable.

\endproclaim
\demo{Proof} $(1) \Rightarrow (2)$, $(2) \Rightarrow (3)$ and $(3) \Rightarrow (4)$ are trivial. $(4) \Rightarrow (5)$ The uniform  topology $\tau_u$ is weaker than the graph topology $\tau_\Gamma$ on $C(X)$. Thus also $(C(X),\tau_u)$ has a countable chain condition. By [MN] $X$ must be compact and metrizable.

$(5) \Rightarrow (1)$ If $X$ is compact and metrizable, the graph  topology $\tau_\Gamma$ on $C(X)$ coincides with the topology $\tau_u$ of uniform convergence, which is second countable if $X$ is compact and metrizable.

\enddemo
\bigskip
\proclaim{Proposition 2.2'} Let $X$ be a Tychonoff space and $Y$ be a metrizable space. The following are equivalent:

(1) $(C(X,Y),\tau_{\Gamma})$ is second countable;

(2) $(C(X,Y),\tau_{\Gamma})$ has a countable network;

(3) $(C(X,Y),\tau_{\Gamma})$ is separable;

(4) $X$ is compact and  metrizable and $Y$ is separable.
\endproclaim

\demo{Proof} $(1)\Rightarrow (2)$ and $(2) \Rightarrow (3)$ are trivial. $(3) \Rightarrow (4)$ The uniform topology $\tau_u$ is weaker than the graph topology $\tau_{\Gamma}$ on $C(X,Y)$. Thus also $(C(X,Y),\tau_u)$ is separable. It is the well known fact that $X$ must be compact and metrizable and $Y$ must be separable.

$(4) \Rightarrow (1)$ If $X$ is compact and metrizable and $Y$ is separable and metrizable, then $(C(X,Y),\tau_u)$ is separable and metrizable. Thus $(C(X,Y),\tau)$ is second countable and $\tau_u = \tau_{\Gamma}$.

\enddemo
\bigskip
\bigskip
\proclaim{Proposition 2.3} Let $X$ be a topological space. Then $(C(X),\tau_\Gamma)$ is a Baire space.

\endproclaim
\demo{Proof} Let $\{G_n: n \in \omega\}$ be a sequence of open dense sets in $(C(X),\tau_\Gamma)$. W.l.o.g we can suppose that $G_{n+1} \subset G_n$ for every $n \in \omega$. We want to prove that $\cap_{n \in \omega} G_n$ is a dense set in $(C(X),\tau_\Gamma)$. Let $U$ be an open set in $(C(X),\tau_\Gamma)$. Define a sequence $\{f_n: n \in \omega\}$ in $C(X)$ and a sequence $\{\epsilon_n: n \in \omega\}$ in  $LSC(X,(0,1))$ such that
\bigskip

\centerline{$B(f_0,\epsilon_0) \subset U \cap G_0$, $B(f_n,\epsilon_n) \subset  B(f_{n-1},\epsilon_{n-1}/2) \cap G_n$ for every $n \ge 1$}
\bigskip

and $\epsilon_n < 1/2^n$ for every $n \in \omega$. There is a continuous function $f$ such that $\{f_n: n \in \omega\}$  $\tau_u$-converges to $f$.
We will show that $f \in U \cap \cap_{n \in \omega} G_n$.

Let $n \in \omega$. We show that $f \in B(f_n,\epsilon_n)$. For every $m > n$ we have $f_m \in B(f_n,\epsilon_n/2)$. Let $x \in X$. The $\tau_u$-convergence of $\{f_n: n \in \omega\}$ to $f$ implies that there is $n_1 > n$ such that $\mid f_m(x) - f(x) \mid < \epsilon_n(x)/2$ for every $m \ge n_1$. We have
\bigskip
\centerline{$\mid f_n(x) - f(x) \mid \le \mid f_n(x) - f_{n_1}(x) \mid + \mid f_{n_1}(x) - f(x) \mid < \epsilon_n(x)/2 + \epsilon_n(x)/2$.}
\bigskip

Thus $f \in B(f_n,\epsilon_n)$ and $B(f_n,\epsilon_n) \subset G_n$ for every $n \in \omega$; i.e. $f \in U \cap \cap_{n \in \omega} G_n$.

\enddemo
\bigskip
\bigskip
\proclaim{Proposition 2.4} Let $X$ be a $T_1$  topological space. Let $Q$ be a subset of $C(X)$. The following are equivalent:

(1) $Q$ is compact in $(C(X),\tau_\Gamma)$;

(2) $Q$ is countably compact in $(C(X),\tau_\Gamma)$;

(3) $Q$ is pseudocompact in $(C(X),\tau_\Gamma)$.

\endproclaim
\demo{Proof} $(1) \Rightarrow (2)$ and $(2) \Rightarrow (3)$ are trivial. We prove $(3) \Rightarrow (1)$. If $Q$ is pseudocompact in $(C(X),\tau_\Gamma)$, the uniform space $(Q,\Cal U_\Gamma)$ is totally bounded [En]. Moreover, $(C(X),\Cal U_\Gamma)$ is complete by Proposition 1.3. If $Q$ is closed in $(C(X),\tau_\Gamma)$, then the uniform space $(Q,\Cal U_\Gamma)$ is complete and totally bounded and hence compact. So it suffices to show that $Q$ is closed in $(C(X),\tau_\Gamma)$.

If $(Q,\tau_\Gamma)$ is pseudocompact, then so is $(Q,\tau_u)$, because any $\tau_u$-continuous function $f: Q \to R$ is $\tau_\Gamma$-continuous and therefore bounded. Thus $(Q,\tau_u)$ is pseudocompact and metrizable and hence compact. So $Q$ is closed in $(C(X),\tau_u)$. Thus $Q$ is closed also in $(C(X),\tau_\Gamma)$. This completes the proof.

\enddemo
\bigskip
\bigskip
\bigskip
\heading 3. Cardinal invariants
\endheading

\bigskip
\bigskip
\proclaim{Theorem 3.1} Let $X$ be a Tychonoff space. Then for every closed discrete set $D$ in $X$ we have $\mid D \mid \le \chi(C(X))$, where $\mid D \mid$ is the cardinality of $D$.

\endproclaim
\demo{Proof}

\enddemo

\enddocument